\documentclass[11pt]{amsart}
\usepackage{tikz-cd}
\usepackage{booktabs}
\usepackage[a4paper,margin=28mm]{geometry}
\usepackage{amsmath,amssymb,mathtools}
\usepackage{microtype,enumitem,xurl}
\usepackage[colorlinks=true,linkcolor=blue,citecolor=blue,urlcolor=blue]{hyperref}
\newtheorem{theorem}{Theorem}[section]
\newtheorem{proposition}[theorem]{Proposition}
\newtheorem{lemma}[theorem]{Lemma}

\newtheorem*{transferthm}{Homotopy transfer}
\theoremstyle{definition}
\newtheorem{definition}[theorem]{Definition}

\theoremstyle{remark}
\newtheorem{remark}[theorem]{Remark}

\numberwithin{equation}{section}

\DeclareMathOperator{\im}{im}
\DeclareMathOperator{\id}{id}
\DeclareMathOperator{\Sh}{Sh}

\newcommand{\Q}{\mathbb Q}

\title{The pentagonal tensor and the fourth obstruction to formality}

\date{November 25, 2026}

\thanks{The second author was supported in part by the National Natural Science Foundation of China (Grant No. 12171165).}
\begin{document}
\author{Chaoshun Yang}
\address{School of Mathematical Sciences,
         Nankai University,
         Tianjin 300071, People's Republic of China}

\email{1296736185@qq.com}

\author{Xiugui Liu}
\address{School of Mathematical Sciences and LPMC,
         Nankai University,
         Tianjin 300071, People's Republic of China}

\email{xgliu@nankai.edu.cn}

\begin{abstract}
Let \(A\) be a connected graded commutative
differential graded algebra over \(\mathbb Q\) whose cohomology
is a finite-dimensional Poincar\'e duality algebra of formal
dimension \(N\ge2\).
Assuming that its Bianchi-Massey tensor vanishes, we prove
that the normalized pentagonal tensor of Nagy and Nordstr\"om, in 
[preprint (2021), arXiv:2105.13660], vanishes if and only if the fourth Harrison obstruction
vanishes after eliminating the third operation.
This gives a fourth-order counterpart to the third-order
comparison of Fiorenza and L\^e, in [Ann. Sc. Norm. Super.
Pisa Cl. Sci. (5), published online 7 March 2025], and answers
the vanishing question in L\^e's proposed comparison,
in [preprint (2026), arXiv:2603.01219], under the hypothesis that the
Bianchi-Massey tensor vanishes.
The result applies to closed, connected, oriented manifolds
with no upper dimension bound.
We also extend Zhou's sphere-bundle formality criterion
in [Algebr. Geom. Topol. 25 (2025), pp. 2297--2315] to two-stage oriented sphere bundles
with positive odd fibre dimensions over a closed, connected,
oriented, rationally formal base, provided that the second
Euler class is pulled back from the original base.
Their total spaces are rationally formal precisely when the
Bianchi-Massey and normalized pentagonal tensors vanish.
Our method develops the model-based strategy in Zhou's proof
 within a commutative bar-cobar
framework, providing a starting point for constructing
higher tensors and comparing them with higher-order
formality obstructions.
\end{abstract}
\maketitle

\section{Introduction}

A commutative differential graded algebra (CDGA) is formal if it
is weakly equivalent to its cohomology algebra \cite{Sullivan,FHT}.
The Bianchi-Massey tensor of Crowley and Nordstr\"om \cite{CN}
organizes quadratic cohomology relations into a choice-independent
obstruction containing information beyond individual triple Massey
products. Zhou used it to characterize formality of orientable
sphere bundles over closed formal manifolds
\cite[Theorem 1.1]{Zhou}.

Nagy and Nordstr\"om introduced the degree-\(-2\) pentagonal tensor,
which is related to fourfold Massey products \cite{NN}.
Under Poincar\'e duality, vanishing of the Bianchi-Massey tensor
permits cochain choices with vanishing uniform triple product;
the pentagonal tensor is independent of choices satisfying this
normalization. They proved that a closed \((n-1)\)-connected
manifold of dimension at most \(5n-2\), where \(n\ge2\), is formal
if and only if its Bianchi-Massey and normalized pentagonal
tensors vanish \cite[Theorem 5.4]{NN}.

The successive obstructions to formality can be expressed through
a transferred minimal \(C_\infty\)-structure on \(H=H^*(A)\).
The first obstruction is the Harrison class of \(m_3\);
after eliminating \(m_3\), the next is the class of \(m_4\)
\cite{Kadeishvili,Le}.
Fiorenza and L\^e established the third-order comparison below
\cite[Theorem 4.4 and Remark 4.5(2)]{FL}, under the Hodge
hypotheses of \cite[Definitions 2.1-2.2]{FL}.

\begin{theorem}[Fiorenza, L\^e]\label{thm:FL}
Let \((A,d)\) be a CDGA over a field \(\Bbbk\) of characteristic
zero, concentrated in degrees \(0,\ldots,N\), whose cohomology
\(H\) is a finite-dimensional Poincar\'e duality algebra of
formal dimension \(N\).
For an orientation \(\varepsilon:H^N\to\Bbbk\), put
\(\langle a,b\rangle=\varepsilon([ab])\) when
\(|a|+|b|=N\), and zero otherwise.
Suppose that there is a degree-\(-1\) map \(h:A\to A\) such that,
writing \(p=1-dh-hd\),
\[
\begin{gathered}
h^2=0,\qquad hdh=h,\qquad dhd=d,\\
\langle\operatorname{im}h,\operatorname{im}h\rangle
=
\langle\operatorname{im}p,\operatorname{im}h\rangle
=0.
\end{gathered}
\]
Then the Bianchi-Massey tensor and the first Harrison obstruction
\(\kappa_3(A)=[m_3]\in H_{\mathrm{Harr}}^{3,-1}(H,H)\)
determine one another. In particular,
\[
\mathcal F_A=0
\quad\Longleftrightarrow\quad
\kappa_3(A)=0.
\]
\end{theorem}

We address the fourth-order comparison posed in
\cite[Section 7(1)]{Le} when the Bianchi-Massey tensor vanishes.
Our result requires only Poincar\'e duality on cohomology and
has no upper dimension bound.
Under these hypotheses, \(\mathcal F_A=0\) also implies
\(\kappa_3(A)=0\), so the third operation can be eliminated.
We prove this implication in Proposition \ref{prop:third}
using Crowley-Nordstr\"om's comparison between the
Bianchi-Massey tensor and uniform triple products.

We now state the main result. A finite-dimensional graded commutative algebra \(H\) is a Poincar\'e duality algebra of formal dimension \(N\) if it has an orientation \(\varepsilon:H^N\to\mathbb Q\) for which \((x,y)\mapsto\varepsilon(xy)\) gives nondegenerate pairings between \(H^i\) and \(H^{N-i}\). We extend \(\varepsilon\) by zero in other degrees. 
\begin{theorem}\label{thm:main}
Let \((A,D_A)\) be a connected graded CDGA over \(\mathbb Q\). Suppose that \(H=H^*(A)\) is a finite-dimensional Poincar\'e duality algebra of formal dimension \(N\ge2\), and that its Bianchi-Massey tensor vanishes:
\[
\mathcal F_A:\mathcal B^d(H)\longrightarrow H^{d-1},
\qquad \mathcal F_A=0.
\]
There exist cochain choices \(c=(\alpha,\gamma)\) with \(\mathcal T_c=0\). All such choices define the same pentagonal tensor
\[
\mathcal P_A:\mathcal D^d(H)\longrightarrow H^{d-2}.
\]
After eliminating the third operation, let
\[
\kappa_4(A)=[m_4]\in H_{\mathrm{Harr}}^{4,-2}(H,H)
\]
be the fourth obstruction. The following conditions are equivalent:
\begin{enumerate}
\item \(\mathcal P_A=0\) in every degree;
\item \(\varepsilon\mathcal P_A|_{\mathcal D^{N+2}(H)}=0\);
\item \(\kappa_4(A)=0\);
\item the minimal \(C_\infty\)-model of \(A\) is isomorphic, by a strictly unital \(C_\infty\)-isomorphism with identity linear component, to a structure in which \(m_3=m_4=0\).
\end{enumerate}
\end{theorem}
The spaces \(\mathcal B(H)\) and \(\mathcal D(H)\), and the tensors \(\mathcal F\), \(\mathcal T\), and \(\mathcal P\), are defined in Section \ref{sec:prelim}. Moreover, we have two applications of Theorem \ref{thm:main}.
\begin{theorem}\label{cor:manifolds}
Let \(M\) be a closed, oriented, connected smooth manifold of dimension \(N\ge2\), and suppose that \(\mathcal F_M=0\). Then the normalized pentagonal tensor is well defined, and
\[
\mathcal P_M=0
\quad\Longleftrightarrow\quad
\varepsilon\mathcal P_M|_{\mathcal D^{N+2}(H^*(M;\mathbb Q))}=0
\quad\Longleftrightarrow\quad
\kappa_4(M)=0.
\]
Equivalently, its minimal rational \(C_\infty\)-model is isomorphic, by a strictly unital \(C_\infty\)-isomorphism with identity linear component, to a structure with \(m_3=m_4=0\).
\end{theorem}

Our second application extends the one-sphere condition in Zhou's theorem to a class of two-stage sphere bundles.

\begin{theorem}\label{cor:spheres}
Let \(B\) be a closed, oriented, connected, rationally formal smooth manifold, and let
\[
S^p\longrightarrow E_1\xrightarrow{\pi_1}B,
\qquad
S^q\longrightarrow X\xrightarrow{\pi_2}E_1
\]
be oriented sphere bundles, where \(p,q\ge1\) are odd. Suppose that the Euler class
\[
e(\pi_2)=\pi_1^*e_2
\quad\text{for some }e_2\in H^{q+1}(B;\mathbb Q).
\]
Put \(N=\dim B+p+q\). The following conditions are equivalent:
\begin{enumerate}
\item \(X\) is rationally formal;
\item \(\mathcal F_X=0\) and the resulting normalized tensor \(\mathcal P_X\) vanishes;
\item \(\mathcal F_X=0\) and \(\varepsilon\mathcal P_X|_{\mathcal D^{N+2}(H^*(X;\mathbb Q))}=0\).
\end{enumerate}
\end{theorem}

The intermediate space \(E_1\) need not be formal, so we cannot use Zhou's theorem directly.

Our proof of Theorem \ref{thm:main} is inspired by Zhou's
construction of a model of cohomology and a map from that model
to the original CDGA \cite[proof of Lemma 3.4]{Zhou}.
We use the commutative bar-cobar model
\(\rho:(\mathcal Q,d)\to(H,0)\)
\cite[Section 4.3]{Kadeishvili}.
A transferred \(C_\infty\)-quasi-isomorphism to \(A\) supplies
an evaluation map \(\Phi:\mathcal Q\to A\) of graded algebras.
Write \(\mathcal Q_{\le r}\) for the subalgebra generated in
bar lengths at most \(r\).
After eliminating \(m_3,\ldots,m_{r-1}\), the map \(\Phi\)
commutes with the differential on \(\mathcal Q_{\le r-1}\).
Its defect \(\mathfrak E=D_A\Phi-\Phi d\) on generators of
length \(r\) represents \(m_r\), up to the suspension sign
(Remark \ref{rem:higher-defect}).
If all higher operations vanish, \(\Phi\) is a CDGA
quasi-isomorphism and exhibits the formality of \(A\).

At arity four, we identify the pentagonal relation space with
the closed elements in an explicit subspace of \(\mathcal Q\).
We construct maps
\[
 Z:\mathcal D(H)\longrightarrow\mathcal Q_{\le3},
 \qquad
 \Gamma:\mathcal D(H)\longrightarrow\mathcal Q_{\le4},
 \qquad d\Gamma=Z,
\]
of degrees \(-2\) and \(-3\), respectively.
Both directions of the theorem follow from the identity
\[
 \Phi Z(z)=D_A\Phi\Gamma(z)-\mathfrak E\Gamma(z).
\]
With \(m_3=0\), the left side represents the normalized
pentagonal tensor. Eliminating \(m_4\) makes the defect term
vanish, proving that the tensor is zero.

For the converse, we use
Poincar\'e duality to detect the Harrison class of \(m_4\).
The key technique is the linearization (Definition \ref{def:linearization}) of the bar-cobar model.
The defect \(\mathfrak E\) vanishes on generators of bar lengths
two and three. Consequently, among the terms of \(\Gamma(z)\),
only those containing a length-four generator contribute to
\(\mathfrak E\Gamma(z)\). Our linearization retains precisely this part.
Thus the proof uses a model of cohomology to obtain the
fourth-order comparison, rather than directly extending
the third-order cochain comparison of Fiorenza-L\^e.

The bar-cobar model suggests a direct approach to constructing
higher tensors and comparing them with higher-order obstructions
to formality. The defect identity applies in every arity and
provides a starting point for extending our method beyond the
fourth obstruction.
We discuss the additional steps required for such extensions
in Section \ref{sec:applications}, together with the finite
obstruction range for models
\[
 A=R\otimes\Lambda(u_1,\ldots,u_s),
 \qquad D_AR=0,\qquad D_Au_i\in R,
\]
where all \(u_i\) have positive odd degree.

The paper is organized as follows.
Section \ref{sec:prelim} recalls the relevant tensors,
the commutative bar-cobar construction, and Harrison obstruction theory.
Section 3 constructs the evaluation map, proving that vanishing of the fourth obstruction implies
vanishing of the pentagonal tensor.
Section \ref{sec:converse} proves the converse using linearization,
Harrison homology, and Poincar\'e duality.
Finally, Section \ref{sec:applications} proves the applications to
manifolds and sphere bundles and discusses possible extensions
to higher obstructions and higher tensors.

\subsection*{AI Usage Statement}

All research objectives and the ideas for constructing the bar-cobar model in this paper were provided by the authors. Some of the proof ideas arose from multiple rounds of interaction between the authors and Chatgpt-6 Astra. The following interactions deserve particular emphasis:

\begin{itemize}
    \item In Section \ref{sec:converse}, inspired by the paper \cite{CN}, we proposed using Poincar\'e duality to detect the obstruction class, but a new functional needed to be constructed based on the relation space. AI helped construct the new functional and proposed using a linearization idea to prove Proposition \ref{prop:converse-surjection} and \ref{prop:converse}. The initial proof constructed by AI used the method of averaging to deal with terms in \(U_2U_3\) (\(U_i\) is defined in Subsection \ref{Def:j,g}), we made use of Proposition \ref{prop:eval-defect} and led to a simplification.

    \item Combining all the results of this paper, the authors summarized the overall method in Section \ref{sec:applications} and proposed an outlook on higher-order problems. Using AI, we constructed a sixth-order tensor according to this method and obtain an fifth-order comparison, thereby verifying the feasibility of this method.
\end{itemize}
We have carefully verified all arguments and take full responsibility for the content of this work.

\section{Preliminaries\label{sec:prelim}}

\subsection{Graded conventions and cochain choices}

Unless explicitly stated otherwise, all vector spaces and tensor products are over \(\mathbb Q\), and all CDGAs are nonnegative and connected, with \(A^0=\mathbb Q1\). We write \(Z^*(A)=\ker D_A\) and \(H=H^*(A)\). 

We adopt the notation of Nagy-Nordstr\"om \cite[Section 4.1]{NN}. The graded symmetric and graded antisymmetric powers are denoted by \(\mathcal G^rH\) and \(\not{\mathcal G}^rH\), respectively. Thus
\[
x\odot y=(-1)^{|x||y|}y\odot x,
\qquad
x\wedge y=-(-1)^{|x||y|}y\wedge x.
\]
The symbol \(\odot\) denotes a formal product in the graded symmetric algebra. In particular, the formal product \(x\odot y\) is not the cohomology class \(xy\). We reserve \(\Lambda\) for a free graded commutative algebra.

Let
\[
\mu:\mathcal G^2H\longrightarrow H,
\qquad \mu(x\odot y)=xy,
\qquad E=\ker\mu.
\]
Let \(K[-]\) be the kernel of full graded symmetrization. More explicitly,
\[
K[H\otimes E]
=\ker\bigl(H\otimes E\longrightarrow\mathcal G^3H\bigr),
\qquad x\otimes e\longmapsto x\odot e.
\]

\begin{definition}\label{def:choices}
A cochain choice is a pair \(c=(\alpha,\gamma)\) consisting of maps
\[
\alpha:H\longrightarrow Z^*(A),\qquad |\alpha|=0,
\qquad
\gamma:E\longrightarrow A,\qquad |\gamma|=-1,
\]
such that \([\alpha(x)]=x\), \(\alpha(1)=1\), and
\[
D_A\gamma=\alpha^2|_E,
\qquad
\alpha^2:\mathcal G^2H\longrightarrow A,
\quad \alpha^2(x\odot y)=\alpha(x)\alpha(y).
\]
\end{definition}

Notice that such choices always exist. We use \(c\) as in \cite[Definition 4.1]{NN}.

\subsection{The uniform triple product and the Bianchi-Massey tensor}

For \(c=(\alpha,\gamma)\), define a degree-\(-1\) map on \(H\otimes E\) by
\[
(\alpha\gamma)(x\otimes e)=(-1)^{|x|}\alpha(x)\gamma(e).
\]
Its differential satisfies
\[
D_A(\alpha\gamma)(x\otimes e)=\alpha^3(x\odot e),
\]
where \(\alpha^3:\mathcal G^3H\to A\) is multiplication of the three chosen representatives. It therefore induces the uniform triple product (called the \emph{uniform Massey triple product} in
\cite[Definition 2.3]{Zhou}) 
\[
\mathcal T_c:K[H\otimes E]^d\longrightarrow H^{d-1},
\qquad
\mathcal T_c(k)=[(\alpha\gamma)(k)].
\]
This is \cite[Definition 4.3]{NN}. In general \(\mathcal T_c\) depends on \(c\); the condition \(\mathcal T_c=0\) will be our normalization.

The Bianchi-Massey domain is
\[
\mathcal B(H)=K[\mathcal G^2E]
=\ker\bigl(\mathcal G^2E\longrightarrow\mathcal G^4H\bigr).
\]
Following the averaged convention of Crowley-Nordstr\"om, used also in \cite[Definition 4.4]{NN}, set
\[
(\gamma\mathbin\cdot\alpha^2)(e\odot f)
=\frac12\bigl(\gamma(e)\alpha^2(f)
+(-1)^{|e||f|}\gamma(f)\alpha^2(e)\bigr).
\]
It takes closed values on \(\mathcal B(H)\) and defines
\[
\mathcal F_A:\mathcal B^d(H)\longrightarrow H^{d-1},
\qquad
\mathcal F_A=[\gamma\mathbin\cdot\alpha^2]|_{\mathcal B(H)}.
\]
This map is independent of \(\alpha\) and \(\gamma\) and is natural under CDGA morphisms \cite[Section 2.2]{CN}. Zhou uses the corresponding unaveraged expression \cite[Definition 2.2]{Zhou}; his tensor is twice the tensor in the present convention, so its vanishing condition is unchanged.

\begin{proposition}\textup{(\cite[Lemma 2.8]{CN})}\label{prop:uniform}
Suppose that \(H=H^*(A)\) is a finite-dimensional Poincar\'e duality algebra of formal dimension \(N\). The top-degree component
\[
\mathcal F_A:\mathcal B^{N+1}(H)\longrightarrow H^N
\]
determines the uniform triple product up to the changes induced by cochain choices, and conversely. In particular, \(\mathcal F_A=0\) if and only if there exists a cochain choice \(c\) with \(\mathcal T_c=0\).
\end{proposition}

\subsection{The pentagonal tensor}

Define the degree-zero map
\[
\begin{aligned}
\mathfrak m:H\otimes\not{\mathcal G}^{2}\mathcal G^2H
&\longrightarrow\mathcal G^3H\otimes\mathcal G^2H,\\
x\otimes(e\wedge f)&\longmapsto
(x\odot e)\otimes f
-(-1)^{|e||f|}(x\odot f)\otimes e.
\end{aligned}
\]
The pentagonal domain of \cite[Definition 4.6]{NN} is
\[
\mathcal D(H)
=\ker\bigl(\mathfrak m|_{H\otimes\not{\mathcal G}^2E}\bigr).
\]
Its definition uses only the graded cohomology algebra. It records an equality between two formal ways of multiplying a cohomology input by one of two quadratic relations.

For a cochain choice \(c\), the induced square of the odd map \(\gamma\) and its product with \(\alpha\) are
\[
\begin{aligned}
\gamma^2:\not{\mathcal G}^2E&\longrightarrow A,
&\gamma^2(e\wedge f)&=(-1)^{|e|}\gamma(e)\gamma(f),\\
\alpha\gamma^2:H\otimes\not{\mathcal G}^2E&\longrightarrow A,
&(\alpha\gamma^2)(x\otimes(e\wedge f))
&=(-1)^{|e|}\alpha(x)\gamma(e)\gamma(f).
\end{aligned}
\]
Both maps have degree \(-2\). The Leibniz rule gives
\[
\begin{aligned}
D_A(\alpha\gamma^2)(x\otimes(e\wedge f))
={}&(-1)^{|x|+|e|}\alpha(x)\alpha^2(e)\gamma(f)\\
&-(-1)^{|x|}\alpha(x)\gamma(e)\alpha^2(f).
\end{aligned}
\]
Equivalently, \(D_A(\alpha\gamma^2)=\mathcal A\mathfrak m\) on \(H\otimes\not{\mathcal G}^2E\), where
\[
\mathcal A:\mathcal G^3H\otimes E\longrightarrow A,
\qquad
\mathcal A(u\otimes f)=(-1)^{|u|}\alpha^3(u)\gamma(f).
\]
Thus \(\alpha\gamma^2\) is closed on \(\mathcal D(H)\).

\begin{definition}\textup{(\cite[Definition 4.8]{NN})}\label{def:pentagonal}
The pentagonal tensor associated to \(c\) is the degree-\(-2\) map
\[
\mathcal P_c:\mathcal D^d(H)\longrightarrow H^{d-2},
\qquad
\mathcal P_c(z)=[(\alpha\gamma^2)(z)].
\]
\end{definition}

\begin{proposition}\textup{(\cite[Corollary 5.2]{NN})}\label{prop:normalized}
If there exists a cochain choice \(c\) with \(\mathcal T_c=0\), then every such choice gives the same map \(\mathcal P_c\). Denote it by \(\mathcal P_A\). This normalized tensor is invariant under CDGA quasi-isomorphisms, with its domain and target identified through the induced cohomology isomorphisms.
\end{proposition}

By Propositions \ref{prop:uniform} and \ref{prop:normalized}, the hypothesis \(\mathcal F_A=0\) in Theorem \ref{thm:main} therefore supplies a well-defined normalized pentagonal tensor. 

\subsection{The commutative bar-cobar construction}
\label{subsec:bar-cobar}

We use the bar-cobar framework of
\cite[Sections 4.1-4.3]{Kadeishvili}, with the following suspension
and sign conventions. For a graded vector space \(V\), let
\(s_V:V\to sV=V[1]\) have degree \(-1\).
For homogeneous \(w_i\in sV\), the shuffle product is
\[
\begin{aligned}
&\mu_{\mathrm{sh}}\bigl(
(w_1\otimes\cdots\otimes w_p)\otimes
(w_{p+1}\otimes\cdots\otimes w_{p+q})\bigr)\\
&\qquad=\sum_{\sigma\in\operatorname{Sh}(p,q)}
\epsilon(\sigma;w)\,
w_{\sigma^{-1}(1)}\otimes\cdots\otimes w_{\sigma^{-1}(p+q)}.
\end{aligned}
\]
Here \(\sigma(i)\) is the new position of the \(i\)-th input,
\(\sigma(1)<\cdots<\sigma(p)\) and
\(\sigma(p+1)<\cdots<\sigma(p+q)\), and
\(\epsilon(\sigma;w)\) is the Koszul sign of this rearrangement.
Put
\[
L^c(sV)=\bigoplus_{r\geq1}L_r^c(sV),
\qquad
L_r^c(sV)=(sV)^{\otimes r}/\operatorname{Sh}_r(sV),
\]
where \(\operatorname{Sh}_r(sV)\) is spanned by shuffle products
of two nonempty tensor words, with Koszul signs computed in \(sV\).
Write \([w_1|\cdots|w_r]\) for a class in this quotient.
The cofree conilpotent Lie coalgebra \(L^c(sV)\) has cobracket
\[
\begin{gathered}
\lambda([w_1|\cdots|w_r])
=\sum_{i=1}^{r-1}
\bigl(u_i\otimes v_i-(-1)^{|u_i||v_i|}v_i\otimes u_i\bigr),\\
u_i=[w_1|\cdots|w_i],\qquad v_i=[w_{i+1}|\cdots|w_r].
\end{gathered}
\]
Its length-one projection
\(\operatorname{pr}_{sV}:L^c(sV)\to sV\)
is the identity on \(L_1^c(sV)=sV\) and zero on all other summands.
For a degree-\(q\) map \(u:L_r^c(sV)\to sV\), its coderivation
extension \(\mathbf u:L^c(sV)\to L^c(sV)\) is given by
\[
\mathbf u([w_1|\cdots|w_n])
=
\sum_{i=0}^{n-r}
(-1)^{q\sum_{j=1}^{i}|w_j|}
[w_1|\cdots|w_i|
u([w_{i+1}|\cdots|w_{i+r}])|
w_{i+r+1}|\cdots|w_n],
\]
where the \(w_j\in sV\) are homogeneous.
For \(n<r\), the value is zero.

For a nonnegative CDGA \((A,D_A,\mu_A)\) with \(A^0=\mathbb Q1\),
define the degree-one coderivation
\(\mathbf B=\mathbf B_1+\mathbf B_2:L^c(sA)\to L^c(sA)\) by
\[
\begin{aligned}
B_1:sA&\longrightarrow sA,
& B_1(s_Aa)&=s_A(D_Aa),\\
B_2:L_2^c(sA)&\longrightarrow sA,
& B_2([s_Aa|s_Ab])&=(-1)^{|a|}s_A(ab).
\end{aligned}
\]
Writing \(\bar A=A^{>0}\), define its reduced commutative
bar construction by
\[
B_{\mathrm{Com}}(A)
=\bigl(L^c(s\bar A),-\mathbf B\bigr),
\]
a differential graded Lie coalgebra.
In our suspension and sign conventions, this is the reduced
version of Kadeishvili's functor \(\Gamma=Q\circ B\),
where \(Q\) takes indecomposables with respect to the
shuffle product \cite[Section 4.3]{Kadeishvili}.

For a conilpotent DG Lie coalgebra
\((C,d_C,\lambda)\), with \(|d_C|=1\) and
\(\lambda:C\to C\otimes C\) of degree zero, let
\(t:C\to tC\) have degree \(1\). Its commutative cobar is
\[
\begin{aligned}
\Omega_{\mathrm{coLie}}(C)&=\bigl(\Lambda(tC),D_1+D_2\bigr),\\
D_1(tc)&=-t(d_Cc),\\
D_2(tc)&=\tfrac12\mu_{\Lambda(tC)}(t\otimes t)\lambda(c).
\end{aligned}
\]
Here \(\Lambda(tC)\) is the free graded commutative algebra,
and both terms extend as degree-one derivations.
These are the commutative cobar formulas of
\cite[Definition 4.10 and Corollary 4.13]{SW}.
This construction corresponds to Kadeishvili's functor
\(\mathcal A\), dual to the Chevalley-Eilenberg functor
\cite[Section 4.3]{Kadeishvili}.

\subsection{Transferred structures and Harrison obstructions}
\label{subsec:obstructions}

We recall the Harrison-cohomological obstruction theory, following \cite[Section 6]{Kadeishvili}.
For a homogeneous map \(f:H^{\otimes r}\to H\), set
\[
\begin{aligned}
\widehat f&=(-1)^{r(r-1)/2}s_Hf(s_H^{-1})^{\otimes r},\\
\widehat f(s_Hx_1\otimes\cdots\otimes s_Hx_r)
&=(-1)^{\sum_{i=1}^r(r-i)|x_i|}
s_Hf(x_1\otimes\cdots\otimes x_r).
\end{aligned}
\]
\begin{transferthm}[\cite{CG}, Theorems 6, 10 and 12]
Let \(H=H^*(A)\). Choose chain maps
\[
  \alpha:(H,0)\longrightarrow(A,D_A),
  \qquad
  \pi:(A,D_A)\longrightarrow(H,0),
\]
and a map \(h:A\to A\) of degree \(-1\) satisfying
\[
  [\alpha(x)]=x,\qquad \alpha(1_H)=1_A,\qquad
  \pi\alpha=\id_H,\qquad
  D_Ah+hD_A=\id_A-\alpha\pi,\qquad
  h\alpha=\pi h=h^2=0.
\]
Then \(H\) carries a strictly unital minimal \(C_\infty\)-structure
with \(m_2\) the cohomology product, together with a strictly unital
\(C_\infty\)-quasi-isomorphism to \(A\) whose linear term is \(\alpha\).

The transferred operations and morphism admit the following
recursive description
\cite[Lemma 7]{CG}.
Written on tensor powers, their components are maps
\[
 F_r:(sH)^{\otimes r}\longrightarrow sA,\qquad |F_r|=0,
 \qquad
 \widehat m_r:(sH)^{\otimes r}\longrightarrow sH,
 \qquad |\widehat m_r|=1.
\]
Starting with \(F_1=s_A\alpha s_H^{-1}\), for \(r\ge2\) set
\begin{equation}\label{eq:transfer-trees}
\begin{aligned}
 F_r
 &=-s_Ahs_A^{-1}
   \sum_{i=1}^{r-1}B_2\circ(F_i\otimes F_{r-i}),\\
 \widehat m_r
 &=s_H\pi s_A^{-1}
   \sum_{i=1}^{r-1}B_2\circ(F_i\otimes F_{r-i}).
\end{aligned}
\end{equation}
All tensor products of maps use the Koszul convention.
For example,
\[
\begin{aligned}
 F_2
 &=-s_Ahs_A^{-1}B_2\circ(F_1\otimes F_1),\\
 \widehat m_3
 &=s_H\pi s_A^{-1}
   \bigl(
      B_2\circ(F_2\otimes F_1)
      +B_2\circ(F_1\otimes F_2)
   \bigr).
\end{aligned}
\]

Expanding the recursion gives a sum over planar rooted binary
trees, with the inputs ordered from left to right.
Each input carries \(F_1\), each binary vertex carries \(B_2\),
each edge joining two binary vertices carries
\(-s_Ahs_A^{-1}\), and the output carries
\(s_H\pi s_A^{-1}\).
Thus every term of \(\widehat m_r\) contains \(r-1\) copies
of \(B_2\) and \(r-2\) copies of \(h\).
\end{transferthm}
Equivalently, the transferred operations are encoded by the
degree-one coderivation
\[
\mathbf b=\sum_{r\geq2}\mathbf b_r:
L^c(sH)\longrightarrow L^c(sH),
\qquad \mathbf b^2=0,
\]
where \(m_r:H^{\otimes r}\to H\) has degree \(2-r\),
\(b_r=\widehat m_r:L_r^c(sH)\to sH\) has degree \(1\),
and \(m_2=\mu\).
It is strictly unital if \(1_H\) is the binary unit and
\(m_r(x_1,\ldots,x_r)=0\) whenever some \(x_i=1_H\), for \(r\geq3\).

The quasi-isomorphism in the transfer theorem is encoded by the
degree-zero Lie coalgebra morphism
\[
\mathbf F:L^c(sH)\longrightarrow L^c(sA),
\qquad \mathbf B\mathbf F=\mathbf F\mathbf b.
\]
Its degree-zero components are
\[
F_r=\left.(\operatorname{pr}_{sA}\circ\mathbf F)
\right|_{L_r^c(sH)}:L_r^c(sH)\longrightarrow sA,
\]
with \(F_1=s_A\alpha s_H^{-1}\), where
\(\alpha:H\to A\) has degree zero,
\(D_A\alpha=0\), and \([\alpha(x)]=x\).
Strict unitality means \(\alpha(1_H)=1_A\) and
\(F_r([s_Hx_1|\cdots|s_Hx_r])=0\) if some \(x_i=1_H\), for \(r\geq2\).
For degree reasons, these maps restrict to the reduced commutative bar
constructions, with source
\(B_{\mathrm{Com}}(H,\{m_r\})=(L^c(sH^{>0}),-\mathbf b)\).

For \(r\geq1\) and \(q\in\mathbb Z\), define
\[
C_{\mathrm{full}}^{r,q}
=\operatorname{Hom}^{q}(L_r^c(sH),sH),
\]
where \(\operatorname{Hom}^{q}\) denotes maps raising degree by \(q\).
Its subspace \(C^{r,q}\) of normalized cochains consists of maps satisfying
\[
u([s_Hx_1|\cdots|s_Hx_r])=0
\quad\text{whenever some }x_i=1_H.
\]
The Harrison differential on the full complex preserves this subspace:
\[
\begin{aligned}
\delta:C^{r,q}&\longrightarrow C^{r+1,q+1},\\
\delta u&=\left.
\bigl(\operatorname{pr}_{sH}\circ[\mathbf b_2,\mathbf u]\bigr)
\right|_{L_{r+1}^c(sH)}.
\end{aligned}
\]
Here \([\mathbf u,\mathbf v]
=\mathbf u\mathbf v-(-1)^{|u||v|}\mathbf v\mathbf u\).
The corresponding unshifted internal degree is \(j=q-r+1\);
we denote Harrison cohomology in arity \(r\) and internal degree
\(j\) by \(H_{\mathrm{Harr}}^{r,j}(H,H)\).
By Poincar\'e duality and \cite[Corollary 7.14]{HLCoh},
the inclusion of normalized Harrison cochains into the full
complex is a quasi-isomorphism.

The arity-four part of \(\mathbf b^2=0\) gives \(\delta b_3=0\).
Its class is
\[
\kappa_3(A)=[b_3]\in H_{\mathrm{Harr}}^{3,-1}(H,H).
\]
If it vanishes, choose \(\xi_2\in C^{2,0}\) with
\(\delta\xi_2=b_3\). Conjugation together with transport of the morphism
\[
\mathbf b'=e^{\boldsymbol{\xi}_2}\mathbf b
e^{-\boldsymbol{\xi}_2},
\qquad
\mathbf F'=\mathbf F e^{-\boldsymbol{\xi}_2}
\]
gives \(b'_3=0\) and \(\mathbf B\mathbf F'=\mathbf F'\mathbf b'\),
and preserves strict unitality.
Each exponential is a finite sum on every fixed tensor length.
After this conjugation, write \(\mathbf b,\mathbf F\)
for \(\mathbf b',\mathbf F'\), respectively. Then \(\delta b_4=0\), and
\[
\kappa_4(A)=[b_4]\in
\frac{\ker(\delta:C^{4,1}\to C^{5,2})}{\delta C^{3,0}}
\cong H_{\mathrm{Harr}}^{4,-2}(H,H).
\]
Its vanishing allows the same conjugation with
\(\xi_3\in C^{3,0}\), \(\delta\xi_3=b_4\), eliminating \(b_4\)
while retaining \(b_3=0\).

\begin{proposition}\label{prop:third}
Let \(A\) be a nonnegative CDGA with \(A^0=\mathbb Q1\), and suppose that \(H=H^*(A)\) is a finite-dimensional Poincar\'e duality algebra. If \(\mathcal F_A=0\), then \(\kappa_3(A)=0\). One may therefore choose a complete strictly unital minimal \(C_\infty\)-structure on \(H\), and an actual strictly unital \(C_\infty\)-quasi-isomorphism to the original \(A\), with \(m_3=0\).
\end{proposition}

\begin{proof}
Let \((H,\{m_r\})\) be the strictly unital minimal
\(C_\infty\)-structure transferred from \(A\), with the chosen
strictly unital quasi-isomorphism \(\mathbf F\) to \(A\).
For the strict CDGA \((H,0,\mu)\), the choices
\(\alpha_H=\operatorname{id}_H\) and \(\gamma_H=0\)
give a zero uniform triple product and a zero ternary operation.
Since \(\mathcal F_A=0\), \cite[Lemma 2.8(ii)]{CN} shows that
the uniform triple product of \(A\) has the same choice class
as this zero choice.

Applying \cite[Lemma 2.12]{CN} with the above choices on
\((H,0,\mu)\) fixed, we obtain an acyclic DGA \(\mathcal E\)
and a minimal \(A_\infty\)-structure
\((H,\{\widetilde m_r\})\) induced by \(A\times\mathcal E\)
such that \(\widetilde m_3=0\).
Thus the class of \(m_3\) in Hochschild cohomology is zero
(see \cite[Section 1.5]{LodayCH} for the definition). 

In characteristic zero, the Harrison cochain complex is a
direct summand of the Hochschild cochain complex
\cite[Section 2.13]{CPRNW}.
Consequently, \([m_3]_{\mathrm{Harr}}=0\).
By the normalization result recalled above, we may choose
\(\xi_2\in C^{2,0}\) such that
\[
\delta\xi_2=b_3=\widehat m_3.
\]
Conjugation and transport of the morphism
\[
\mathbf b'
=
e^{\boldsymbol{\xi}_2}\mathbf b
e^{-\boldsymbol{\xi}_2},
\qquad
\mathbf F'
=
\mathbf F e^{-\boldsymbol{\xi}_2}
\]
eliminates \(b_3\) and preserves strict unitality and the
morphism equation to \(A\).
\end{proof}

This establishes the low-order input under the cohomological duality hypothesis of Theorem \ref{thm:main}. The proof does not assume a Hodge homotopy of the kind in Theorem \ref{thm:FL}. The subsequent construction uses the ordinary transferred structure and the nondegenerate pairing on \(H\).

\section{From the fourth obstruction to the pentagonal tensor}

We prove the implication \(\kappa_4(A)=0\Rightarrow\mathcal P_A=0\)
in Theorem \ref{thm:main}. We first compute the defect of evaluation
on the free model of \(H\), then construct a cochain in that model
whose differential represents the pentagonal tensor.

\subsection{Evaluation and the fourth operation}\label{Def:j,g}
Apply the construction of
Section \ref{subsec:bar-cobar} to the strict algebra \((H,0,\mu)\).
With \(\bar H=H^{>0}\), this gives
\[
  (\mathcal Q,d)=\Omega_{\mathrm{coLie}}B_{\mathrm{Com}}(H,0,\mu)
  =\bigl(\Lambda(\textstyle\bigoplus_{r\ge1}U_r),D_1+D_2\bigr),
  \qquad U_r=tL_r^c(s_H\bar H).
\]
Write \(\mathcal Q_{\le r}=\Lambda(U_1\oplus\cdots\oplus U_r)\).
The counit is the CDGA map determined on generators by
\begin{equation}\label{eq:cobar-counit}
 \rho:(\mathcal Q,d)\longrightarrow(H,0),\qquad
 \rho(t[s_Hx])=x\quad(x\in\bar H),\qquad
 \rho|_{U_r}=0\quad(r\ge2).
\end{equation}
The counit \(\rho\) is a quasi-isomorphism by the
commutative bar-cobar resolution theorem
\cite[Section 4.3]{Kadeishvili}
(see also \cite[Corollary 11.3.5]{LV}).
Its linear section \(j:H\to\mathcal Q\) is given by \(j(1)=1\)
and \(j(x)=t[s_Hx]\) for \(x\in\bar H\). Thus \(\rho j=\id_H\).

For \(q\in\mathcal G^2\bar H\) put \(e(q)=q-1\odot\mu(q)\).
Every element of \(E\) has this form uniquely, as one sees by
decomposing \(\mathcal G^2H\) according to the number of unit factors.
Define \(g:E\to U_2\) of degree \(-1\) by
\[
 g(e(x\odot y))=(-1)^{|x|+1}t[s_Hx|s_Hy].
\]
The binary shuffle makes this expression graded symmetric in \(x,y\).
The two cobar terms give \(dg(e(q))=j^2(q)-j(\mu(q))=j^2(e(q))\),
where \(j^2(x\odot y)=j(x)j(y)\).

\begin{proposition}[Defect of evaluation]\label{prop:eval-defect}
Take the morphism \(\mathbf F\) of Proposition \ref{prop:third}, with
\(b_3=0\). Let \(\Phi:\mathcal Q\to A\) be the unital graded algebra
map defined by
\[
 \Phi(tu)=s_A^{-1}F_r(u),\qquad u\in L_r^c(s_H\bar H).
\]
The degree-one map \(\mathfrak E=D_A\Phi-\Phi d:\mathcal Q\to A\) satisfies
\begin{equation}\label{eq:forward-defect}
 \mathfrak E|_{\mathcal Q_{\le3}}=0,\qquad
 \mathfrak E(tu)=s_A^{-1}F_1 b_4(u)
 \quad\bigl(u\in L_4^c(s_H\bar H)\bigr).
\end{equation}
Moreover, we have
\(\mathfrak E(tu)=(-1)^{3|x_1|+2|x_2|+|x_3|}\alpha(m_4(x_1,x_2,x_3,x_4))\) for \(u=[s_Hx_1|\cdots|s_Hx_4]\).
In particular \(\Phi|_{\mathcal Q_{\le4}}\) is a CDGA map when \(b_4=0\).
\end{proposition}
\begin{proof}
The \(F_r\) are defined on the shuffle quotients, so \(\Phi\) is well defined. Since the bar differential
of \((H,0,\mu)\) is \(-\mathbf b_2\), its cobar differential has
\(D_1(tu)=t(\mathbf b_2u)\). And we have that \(D_A(\Phi(tu))=s_A^{-1}B_1F_r(u)\). Choose a word representing \(u=[w_1|\cdots|w_r]\)
and
write \(u_i=[w_1|\cdots|w_i]\) and \(v_i=[w_{i+1}|\cdots|w_r]\).
The two terms of the antisymmetrized cobracket give equal contributions
to \(D_2\), cancelling its factor \(1/2\). Thus
\[
 D_2(tu)=\sum_{i=1}^{r-1}(-1)^{|u_i|}(tu_i)(tv_i).
\]
Put \(a_i=s_A^{-1}F_i(u_i)\) and \(a'_i=s_A^{-1}F_{r-i}(v_i)\).
The target multiplication has the opposite sign in this term, since
\[
 s_A^{-1} B_2([F_i(u_i)|F_{r-i}(v_i)])
    =(-1)^{|a_i|}a_i a'_i=(-1)^{|u_i|+1}a_i a'_i.
\]
It follows that
\[
 \mathfrak E(tu)=s_A^{-1}\big(B_1F_r(u)+\sum_{i=1}^{r-1} B_2([F_i(u_i)|F_{r-i}(v_i)])-F_{r-1}(\mathbf b_2u)\big),
\]
with the last term omitted for \(r=1\). Projecting
\(\mathbf B\mathbf F=\mathbf F\mathbf b\) to \(sA\) on this word
therefore gives the identity
\begin{equation}\label{eq:all-defects}
 \mathfrak E(tu)
   =s_A^{-1}\sum_{k=3}^{r}F_{r-k+1}(\mathbf b_k u).
\end{equation}
For \(r\le4\), only \(F_1b_4\) can remain. The sign in the statement
follows from \(b_4=\widehat m_4\) and \(F_1=s_A\alpha s_H^{-1}\).
Finally, for homogeneous \(p,q\in\mathcal Q\),
\[
 \mathfrak E(pq)=\mathfrak E(p)\Phi(q)+(-1)^{|p|}\Phi(p)\mathfrak E(q).
\]
This extends the computation on generators to their products.
\end{proof}

\begin{remark}\label{rem:higher-defect}
For \(n\ge3\), if \(b_3,\ldots,b_{n-1}\) have been eliminated, \eqref{eq:all-defects}
gives \(\mathfrak E|_{\mathcal Q_{\le n-1}}=0\) and
\(\mathfrak E(tu)=s_A^{-1}F_1b_n(u)\) for \(u\in L_n^c(s_H\bar H)\).
\end{remark}

\subsection{Exactness in each weight}
For \(t[s_Hx_1|\cdots|s_Hx_r]\in U_r\), define the internal degree
\(p=\sum_i|x_i|\), the weight \(n=r\) and the resolution degree \(a=r-1\).
All three are additive on products, and the cohomological degree is
\(p-a\). Write \(\mathcal Q_a^{(n)}\) for the part of weight \(n\) and
resolution degree \(a\), and let 
\[
    F_n\mathcal Q=\bigoplus_{m\le n}\mathcal{Q}^{(m)},\quad \mathcal{Q}^{(m)}=\bigoplus_j \mathcal{Q}^{(m)}_j.
\]
The two parts of \(d\) preserve the internal degree
and have the following degrees:
\[
 D_1:\mathcal Q_a^{(n)}\longrightarrow\mathcal Q_{a-1}^{(n-1)},
 \qquad
 D_2:\mathcal Q_a^{(n)}\longrightarrow\mathcal Q_{a-1}^{(n)}.
\]
Thus the differential on the weight quotient
\(F_n\mathcal Q/F_{n-1}\mathcal Q=\mathcal Q^{(n)}\) is \(D_2\).

\begin{lemma}\label{lem:weight-exact}
For every \(n\ge2\), the complex \((\mathcal Q^{(n)},D_2)\) is exact
in every internal degree. In particular, indexing homology by
resolution degree, we have
\begin{equation}\label{eq:three-weights}
 H_2(\mathcal Q^{(5)},D_2)=H_2(\mathcal Q^{(4)},D_2)=0,
 \qquad D_2:U_3=\mathcal Q^{(3)}_2\longrightarrow
              \mathcal Q^{(3)}_1\ \text{is injective}.
\end{equation}
\end{lemma}
\begin{proof}
Put \(W=s_H\bar H\) and \(V=W^\vee\), where \(W^\vee\) denotes the dual space of \(W\). Two dualities enter the argument.

First, notice that the following two Hopf algebras 
\[
    (T^C(W),\mu_{sh},\Delta_{dec}),\quad (T(V),\mu_{tensor},\Delta_{H})
\]
are dual in each tensor length under the pairing
\[
\langle w_1\cdots w_n,a_1\cdots a_n\rangle
=\prod_{i=1}^n a_i(w_i).
\]
More explicitly, we have two identities:

\begin{equation}\label{dual1}
    \langle c,ab\rangle=\langle\Delta_{\mathrm{dec}}c,a\otimes b\rangle 
,\text{ where } a,b\in TV,
\end{equation}
\begin{equation}\label{dual2}
    \begin{aligned}
\langle \mu_{sh}(c_1\otimes c_2),h\rangle 
 &=\sum_{\sigma\in\operatorname{Sh}(p,q)}
      \langle\sigma(c_1c_2),h\rangle\\
 &=\sum_{\sigma\in\operatorname{Sh}(p,q)}
      \langle c_1c_2,\sigma^{-1}h\rangle
  =\langle c_1\otimes c_2,\Delta_Hh\rangle ,\quad h\in TV.
\end{aligned}
\end{equation}

For the notation in (\ref{dual2}), see \cite[Proposition 1.3.4]{LV}. On \(T^{C}(W)\), the deconcatenation coproduct is given by
\[
\begin{aligned}
\Delta_{\mathrm{dec}}(w_1\cdots w_n)
&=\sum_{i=0}^{n}
  (w_1\cdots w_i)\otimes(w_{i+1}\cdots w_n).
\end{aligned}
\]

On \(T(V)\), the product is concatenation and the coproduct is
the algebra map determined by declaring the elements of \(V\) primitive:
\[
\begin{aligned}
\mu_{\mathrm{tensor}}(u\otimes v)&=uv,\\
\Delta_H(a_1\cdots a_n)
&=\prod_{i=1}^{n}(a_i\otimes1+1\otimes a_i),
\qquad a_i\in V.
\end{aligned}
\]

Let \(\mathbb L(V)\) be the free graded Lie algebra on \(V\),
with \(\mathbb L_n(V)\) its part of bracket length \(n\). The evaluation
pairing defined above descends to
\begin{equation}\label{eq:lie-word-dual}
 (L_n^c(W))^\vee\cong\mathbb L_n(V),\qquad
 \langle\lambda u,a\otimes b\rangle=\langle u,[a,b]\rangle ,
\end{equation}
where \(\lambda\) is the cobracket on Lie coalgebra \(L_n^c(W)\). Indeed, let \(T^C_+(W)\) be the subspace of positive tensor length and let
\(\operatorname{Sh}(W)\) be the span of \(\mu_{sh}(c_1\otimes c_2)\) with
both \(c_1\) and \(c_2\) of positive length. Write
\[
\pi:T^C_+(W)\longrightarrow
L^c(W):=T^C_+(W)/\operatorname{Sh}(W)
\]
for the quotient map. By (\ref{dual2}) and the nondegeneracy of the pairing, for
\(h\in V^{\otimes n}\) with \(n>0\),
\[
\langle\operatorname{Sh}(W),h\rangle=0
\quad\Longleftrightarrow\quad
\bar\Delta_Hh=0, \text{ where }\bar\Delta_H\text{ is the reduced coproduct.}
\]
Over \(\mathbb Q\), the primitive elements of \(T(V)\) form the
free graded Lie algebra \(\mathbb L(V)\)
\cite[Corollary 1.3.7]{LV}. Thus evaluation induces the isomorphism
\[
\mathbb L_n(V)\xrightarrow{\ \cong\ }(L_n^c(W))^\vee,
\qquad
h\longmapsto\bigl(\pi c\longmapsto\langle c,h\rangle\bigr),
\]
where \(\mathbb L_n(V)\) denotes the part of bracket length \(n\).

Finally, let \(\tau(x\otimes y)=(-1)^{|x||y|}y\otimes x\).
The Lie cobracket
\(\lambda:L^c(W)\to L^c(W)\otimes L^c(W)\) is induced by
antisymmetrized deconcatenation \cite[\S1.2.4]{BFMT}:
\[
\lambda\pi=(\pi\otimes\pi)(1-\tau)\bar\Delta_{\mathrm{dec}}.
\]
Thus by (\ref{dual1}), we have
\begin{equation}\label{eq:lie-word-dual2}
\begin{aligned}
\langle\lambda u,a\otimes b\rangle
 &=\langle\bar\Delta_{\mathrm{dec}}c,a\otimes b\rangle
   -(-1)^{|a||b|}
    \langle\bar\Delta_{\mathrm{dec}}c,b\otimes a\rangle\\
 &=\langle c,ab-(-1)^{|a||b|}ba\rangle
  =\langle u,[a,b]\rangle .
\end{aligned}
\end{equation}

Second, put \(C=L^c(W)\) and \(\mathfrak g=\mathbb L(V)\),
paired as in \eqref{eq:lie-word-dual}.
Regard \(V\) and \(\mathfrak g\) as homologically graded,
with \(V_i=(W^i)^\vee\).
In the weight quotient only \(D_2\) survives, so we equip
\(C\) and \(\mathfrak g\) with zero internal differential.

Recall that \(t:C\to tC\) raises cochain degree by one.
Write
\[
 s_{\mathfrak g}:\mathfrak g\longrightarrow s_{\mathfrak g}\mathfrak g,
 \qquad |s_{\mathfrak g}a|=|a|+1,
\]
for the suspension in chain degree.
The commutative cobar algebra of \(C\) is
\((\Lambda(tC),D_2)\).
The Chevalley-Eilenberg chain coalgebra of \(\mathfrak g\)
has length-\(k\) part
\[
 C_k^{\mathrm{CE}}(\mathfrak g)
 =\Lambda^k(s_{\mathfrak g}\mathfrak g),
\]
with coproduct determined by declaring
\(s_{\mathfrak g}\mathfrak g\) primitive.

Set \(\langle tu,s_{\mathfrak g}a\rangle=\langle u,a\rangle\).
Pair unequal lengths by zero and equal lengths by
\[
 \left\langle
 (tu_1)\cdots(tu_k),
 (s_{\mathfrak g}a_1)\cdots(s_{\mathfrak g}a_k)
 \right\rangle
 =
 \sum_{\sigma\in S_k}
 \varepsilon(\sigma;s_{\mathfrak g}a_1,\ldots,s_{\mathfrak g}a_k)
 \prod_{i=1}^k\langle u_i,a_{\sigma(i)}\rangle,
\]
where \(\varepsilon(\sigma;s_{\mathfrak g}a_1,\ldots,
s_{\mathfrak g}a_k)\) is the Koszul sign of reordering
\((s_{\mathfrak g}a_1,\ldots,s_{\mathfrak g}a_k)\) as
\((s_{\mathfrak g}a_{\sigma(1)},\ldots,
s_{\mathfrak g}a_{\sigma(k)})\).

This pairing gives the following diagram, with duals taken
weight by weight; compare
\cite[Theorem 4.17 and its proof]{SW}.
\[
\begin{tikzcd}[column sep=large,row sep=large]
 (C,0)
   \arrow[r,mapsto,"{(-)^\vee}"]
   \arrow[d,mapsto,"{(-)^\vee\circ\Omega_{\mathrm{coLie}}}"']
 & (\mathfrak g,0)
   \arrow[d,mapsto,"{C_*^{\mathrm{CE}}}"] \\
 (\Lambda(tC),D_2)^\vee
   \arrow[r,"\cong"']
 & (\Lambda(s_{\mathfrak g}\mathfrak g),
       \partial_{\mathrm{br}}).
\end{tikzcd}
\]

To verify compatibility with the differentials, write
\(\lambda u=\sum u'\otimes u''\).
The defining formulas are
\[
 D_2(tu)=\tfrac12\sum(-1)^{|u'|}(tu')(tu''),
 \qquad
 \partial_{\mathrm{br}}
 \bigl((s_{\mathfrak g}a)(s_{\mathfrak g}b)\bigr)
 =(-1)^{|a|}s_{\mathfrak g}[a,b].
\]
The length-two pairing and graded antisymmetry of \(\lambda\)
give
\[
\begin{aligned}
 \left\langle D_2(tu),
       (s_{\mathfrak g}a)(s_{\mathfrak g}b)\right\rangle
 &=(-1)^{|a|}\langle\lambda u,a\otimes b\rangle\\
 &=\left\langle tu,
       \partial_{\mathrm{br}}
       \bigl((s_{\mathfrak g}a)(s_{\mathfrak g}b)\bigr)
   \right\rangle.
\end{aligned}
\]
Since \(D_2\) extends as a derivation and
\(\partial_{\mathrm{br}}\) as a coderivation,
this identity extends to all lengths.

A product of \(k\) factors of total weight \(n\) has
resolution degree \(n-k\). It follows that
\begin{equation}\label{eq:ce-weight-dual}
 (\mathcal Q_{n-k}^{(n)})^\vee
 \cong C_k^{\mathrm{CE}}(\mathfrak g)^{(n)},
 \qquad D_2^\vee=\partial_{\mathrm{br}}.
\end{equation}
In internal degree \(p\), both sides have chain degree
\(p-n+k\).
The pairing is nondegenerate in each weight.

Finally, put \(R=U(\mathfrak g)=T(V)\). The Chevalley-Eilenberg
resolution identifies \(H_k^{\mathrm{CE}}(\mathfrak g;\mathbb Q)\)
with \(\operatorname{Tor}^R_k(\mathbb Q,\mathbb Q)\), preserving weight
\cite[Propositions 22.3-22.4]{FHT}. The trivial right \(R\)-module
has the weight-preserving free resolution
\[
 0\longrightarrow V\otimes R
 \xrightarrow{\ v\otimes r\,\mapsto\,vr\ }R
 \xrightarrow{\ \varepsilon\ }\mathbb Q\longrightarrow0 ,
\]
since each nonempty word has a unique first letter.
By \eqref{eq:ce-weight-dual},
\((\mathcal Q^{(n)},D_2)\) is exact in every internal degree for
\(n\ge2\). The injectivity in \eqref{eq:three-weights} follows since
\(\mathcal Q_3^{(3)}=0\).
\end{proof}

\subsection{The lift of the pentagonal relations}
Put \(K=K[H\otimes E]\). With the notation of Section \ref{sec:prelim} and Subsection \ref{Def:j,g},
write \(S=(jg)|_K\) and \(Z=(jg^2)|_{\mathcal D(H)}\), where
\[
 (jg)(x\otimes e)=(-1)^{|x|}j(x)g(e),\qquad
 (jg^2)(x\otimes(e\wedge f))=(-1)^{|e|}j(x)g(e)g(f).
\]

\begin{lemma}[Finite lifts of the relation spaces]\label{lem:finite-lift}
There are linear maps \(Y,\Gamma\) of degrees \(-2,-3\), respectively,
with
\begin{equation}\label{eq:forward-lifts}
\begin{aligned}
 Y &:K\longrightarrow U_3, & dY&=S,\\
 \Gamma &:\mathcal D(H)\longrightarrow U_1U_4\oplus U_2U_3\oplus U_4,
 &d\Gamma&=Z .
\end{aligned}
\end{equation}
Any two such choices of \(\Gamma\) differ by \(dT\) for a linear
map \(T:\mathcal D(H)\to U_5\) of degree \(-4\).
\end{lemma}
\begin{proof}
Since \(dj=0\) and \(dg=j^2|_E\), the two Leibniz identities of
Section \ref{sec:prelim} give \(dS=0\) and \(dZ=0\). Explicitly,
\(d(jg)(x\otimes e)=j^3(x\odot e)\) and
\(d(jg^2)=\mathcal A_j\mathfrak m\), where
\(\mathcal A_j:\mathcal G^3H\otimes E\to\mathcal Q\) sends
\(v\otimes f\) to \((-1)^{|v|}j^3(v)g(f)\).
Both representatives map to zero under the counit \(\rho\) of \eqref{eq:cobar-counit}.

For \(n\ge2\), the quotient complex
\[
 F_n\mathcal Q/F_{n-1}\mathcal Q
 \cong (\mathcal Q^{(n)},D_2)
\]
is exact by Lemma \ref{lem:weight-exact}. Thus the inclusion
\(F_1\mathcal Q\hookrightarrow F_n\mathcal Q\) induces an
isomorphism on cohomology. Since \(F_1\mathcal Q=j(H)\) and
\(\rho j=\id_H\), its inverse, after identifying \(F_1\mathcal Q\)
with \(H\) via \(j\), is
\[
 \rho_*:H^*(F_n\mathcal Q,d)\xrightarrow{\;\cong\;}H.
\]

For \(k\in K\), the cocycle \(S(k)\in F_3\mathcal Q\) has
resolution degree one. Since \(\rho S(k)=0\), its class in
\(H^*(F_3\mathcal Q,d)\) is zero. As \(d\) lowers resolution
degree by one, we may choose a preimage \(Y(k)\) of resolution
degree two. The part of \(F_3\mathcal Q\) in that degree is
\[
 \bigoplus_{m\le3}\mathcal Q_2^{(m)}
 =\mathcal Q_2^{(3)}=U_3.
\]
Thus \(Y(k)\in U_3\) and \(dY(k)=S(k)\).

For \(z\in\mathcal D(H)\), decompose \(Z(z)\) by weight as
\[
 Z(z)=Z_5+Z_4,\qquad Z_5\in U_1\mathcal G^2U_2,\quad
 Z_4\in\mathcal G^2U_2 .
\]
The term \(Z_4\) comes from the scalar part of \(j(x)\).
Both terms have resolution degree two. Taking the weight five, four
and three components of \(dZ(z)=0\) gives, respectively,
\[
 D_2Z_5=0,\qquad D_2Z_4=-D_1Z_5,\qquad D_1Z_4=0 .
\]
By \eqref{eq:three-weights} choose \(Y_5\in U_1U_4\oplus U_2U_3\)
with \(D_2Y_5=Z_5\). Then
\[
 D_2(Z_4-D_1Y_5)=-D_1Z_5+D_1D_2Y_5=0,
\]
so the weight four assertion gives \(Y_4\in U_4\) with
\(D_2Y_4=Z_4-D_1Y_5\). Finally \(D_1Y_4\in U_3\) and
\(D_2D_1Y_4=-D_1Z_4+D_1^2Y_5=0\). Injectivity on \(U_3\) gives
\(D_1Y_4=0\). Thus \(\Gamma(z)=Y_5+Y_4\) satisfies \(d\Gamma(z)=Z(z)\).
Choosing the lifts on homogeneous bases gives the stated linear maps.

For two such choices \(\Gamma,\Gamma'\) and
\(z\in\mathcal D(H)\), the difference
\((\Gamma-\Gamma')(z)\) is a cocycle of resolution degree three
in \(F_5\mathcal Q\). Its image under \(\rho\) is zero by
\eqref{eq:cobar-counit}, so it is exact in \(F_5\mathcal Q\).
We may choose a preimage \(T(z)\) of resolution degree four.
The part of \(F_5\mathcal Q\) in that degree is
\[
 \bigoplus_{m\le5}\mathcal Q_4^{(m)}
 =\mathcal Q_4^{(5)}=U_5.
\]
Choosing these preimages on a homogeneous basis of \(\mathcal D(H)\)
gives a linear map \(T:\mathcal D(H)\to U_5\) of degree \(-4\)
with \(\Gamma-\Gamma'=dT\).
\end{proof}

\begin{proposition}\label{prop:forward}
Under the hypotheses of Theorem \ref{thm:main}, if \(\kappa_4(A)=0\),
then \(\mathcal P_A:\mathcal D^d(H)\to H^{d-2}\) vanishes for every \(d\).
\end{proposition}
\begin{proof}
Set \(\alpha=\Phi j\) and \(\gamma=\Phi g\). Then \(
\alpha(x)=s_A^{-1}F_1(s_Hx)\), and \(D_A\gamma=\Phi dg=\alpha^2|_E\).
Thus \(c=(\alpha,\gamma)\) is a cochain choice. Since \(Y\) takes values
in \(U_3\), equations \eqref{eq:forward-lifts} and \eqref{eq:forward-defect}
give \(D_A\Phi Y=(\alpha\gamma)|_K\). Hence \(\mathcal T_c=0\).

If \(\kappa_4(A)=0\), apply the conjugation described in
Section \ref{subsec:obstructions} with \(\delta\xi_3=b_4\).
This leaves \(F_1,F_2\), and hence \(c\), unchanged. For the transported
morphism, still denoted by \(\mathbf F\), we have \(b_3=b_4=0\).
By \eqref{eq:forward-defect}, its evaluation \(\Phi\) commutes with
the differential on \(\mathcal Q_{\le4}\). Consequently
\[
 D_A\bigl(\Phi\Gamma(z)\bigr)=\Phi Z(z)=(\alpha\gamma^2)(z)
 \qquad (z\in\mathcal D(H)).
\]
Taking cohomology gives \(\mathcal P_c=0\) by Definition \ref{def:pentagonal},
and \(\mathcal P_A=\mathcal P_c\) by Proposition \ref{prop:normalized}.
\end{proof}

\section{From the pentagonal tensor to the fourth obstruction}
\label{sec:converse}

Choose the structure and morphism of Proposition \ref{prop:third} with
\(b_3=0\), retaining \(b_4=\widehat m_4\). We use the same
\(\Phi,\alpha,\gamma\) as in the proof of Proposition \ref{prop:forward}. The first
identity in \eqref{eq:forward-lifts} gives
\(D_A\Phi Y=(\alpha\gamma)|_K\), since \(\Phi\) commutes with the
differential on \(\mathcal Q_{\le3}\). Thus \(\mathcal T_c=0\) for
\(c=(\alpha,\gamma)\), and the class of \(\Phi Z(z)\) represents \(\mathcal P_A(z)\).

\subsection{Harrison chains and Poincar\'e duality}

Write \(C_r(H)=H\otimes(\bar H^{\otimes r}/\Sh_r)\) for the
normalized Harrison chains with coefficients in \(H\)
\cite[\S4.2.10]{LodayCH}. Here \(\Sh_r\) is spanned by shuffles of
two nonempty words with sign \(\operatorname{sgn}(\pi)\epsilon(\pi;\mathbf x)\),
where \(\epsilon(\pi;\mathbf x)\) is the Koszul sign on the unshifted
inputs. A chain \(v[x_1|\cdots|x_r]\) has coefficient \(v\in H\)
 and inputs \(x_i\in\bar H\). Write \(C_r(H)^n\) for
the part of total internal degree \(n=|v|+\sum_i|x_i|\).
The boundary \(\partial_r:C_r(H)^n\to C_{r-1}(H)^n\) is
\begin{align*}
 \partial_r\bigl(v[x_1|\cdots|x_r]\bigr)
 &=vx_1[x_2|\cdots|x_r]\\
 &\quad+\sum_{i=1}^{r-1}(-1)^i
 v[x_1|\cdots|x_ix_{i+1}|\cdots|x_r]\\
 &\quad+(-1)^{r+|x_r|(|v|+\sum_{i<r}|x_i|)}
 x_rv[x_1|\cdots|x_{r-1}].
\end{align*}

Recall from Section 1 that \(\varepsilon:H^N\to\mathbb Q\) is the
orientation of \(H\), extended by zero in all other degrees, and that
\((v,y)\mapsto\varepsilon(vy)\) pairs \(H^i\) nondegenerately with
\(H^{N-i}\).

\begin{proposition}\label{prop:converse-pd}
For each \(r\ge1\), evaluation against the orientation induces an
isomorphism
\[
 C^{r,r-3}\xrightarrow{\ \cong\ }(C_r(H)^{N+2})^\vee,\qquad
 \widehat f\longmapsto\lambda_f,\qquad
 \lambda_f(v[x_1|\cdots|x_r])
   =\varepsilon\bigl(vf(x_1,\ldots,x_r)\bigr).
\]
If \(\widehat{f'}=\delta\widehat f\), then
\begin{equation}\label{eq:converse-pd}
 \lambda_{f'}=(-1)^{r+1}\lambda_f\circ\partial_{r+1}.
\end{equation}
\end{proposition}

\begin{proof}
By the suspension convention of Subsection \ref{subsec:obstructions},
\(\widehat f\in C^{r,r-3}\) corresponds to a normalized Harrison map
\(f:H^{\otimes r}\to H\) of internal degree \(-2\).
Let \(L\in(C_r(H)^{N+2})^\vee\). Fix homogeneous inputs
\(x_1,\ldots,x_r\in\bar H\) and put \(s=\sum_i|x_i|\).
Then the map
\[
 H^{N+2-s}\longrightarrow\mathbb Q,\qquad
 v\longmapsto L(v[x_1|\cdots|x_r])
\]
is linear. By Poincar\'e duality, there is a unique element
\(f_L(x_1,\ldots,x_r)\in H^{s-2}\) such that
\[
 L(v[x_1|\cdots|x_r])
 =\varepsilon\bigl(vf_L(x_1,\ldots,x_r)\bigr)
 \qquad\text{for every }v\in H^{N+2-s}.
\] The shuffle relations give the
Harrison condition, and extension by zero on unit inputs gives
normalization. This constructs the inverse.

Desuspending \([\mathbf b_2,\boldsymbol{\widehat f}]\) gives
\begin{align*}
 f'(x_1,\ldots,x_{r+1})
 &=(-1)^{r+1}x_1 f(x_2,\ldots,x_{r+1})\\
 &\quad+\sum_{i=1}^r(-1)^{r+i+1}
       f(x_1,\ldots,x_ix_{i+1},\ldots,x_{r+1})\\
 &\quad+f(x_1,\ldots,x_r)x_{r+1}.
\end{align*}
Multiply by \(v\) and apply \(\varepsilon\). In the last term, moving
\(x_{r+1}\) to the left gives
\((-1)^{|x_{r+1}|(|v|+\sum_{i=1}^r|x_i|)}\), since \(|f|=-2\).
Comparison with \(\partial_{r+1}\) gives \eqref{eq:converse-pd}.
\end{proof}

Since \(\delta b_4=0\), equation \eqref{eq:converse-pd} gives
\(\lambda_{m_4}\partial_5=0\). Denote the induced functional on
\(H_4(C_\bullet(H)^{N+2})\) by \(\overline\lambda_{m_4}\).

\subsection{From pentagonal relations to Harrison homology}

\begin{proposition}\label{prop:converse-kernel}
The degree-\(-2\) map
\[
 jg^2:H\otimes\not{\mathcal G}^2E\longrightarrow j(H)\mathcal G^2U_2
\]
is an isomorphism. Its restriction \(Z\) to \(\mathcal D(H)\)
identifies the pentagonal relations with the closed elements:
\begin{equation}\label{eq:converse-kernel}
 Z(\mathcal D(H))=j(H)\mathcal G^2U_2\cap\ker d .
\end{equation}
\end{proposition}

\begin{proof}
Recall that \(j:H\to j(H)=\mathbb Q1\oplus U_1\) and \(g:E\to U_2\)
are isomorphisms of degrees \(0\) and \(-1\), respectively.
The sign in \(g^2(e\wedge f)=(-1)^{|e|}g(e)g(f)\) identifies
\(\not{\mathcal G}^2E\) with \(\mathcal G^2U_2\) up to the degree
shift. Since \(\mathcal Q\) is freely generated by the \(U_i\),
multiplication identifies \(j(H)\otimes\mathcal G^2U_2\) with
\(j(H)\mathcal G^2U_2\). This proves the first assertion.

By the proof of Lemma \ref{lem:finite-lift},
\(d(jg^2)=\mathcal A_j\mathfrak m\), where
\(\mathcal A_j(u\otimes f)=(-1)^{|u|}j^3(u)g(f)\).
For \(0\le i\le3\), \(j^3\) maps
\(1^{\odot(3-i)}\odot\mathcal G^i\bar H\) isomorphically onto
\(\mathcal G^iU_1\). Moreover, \(g(f)\in U_2\), so
\(\mathcal A_j\) is injective. It follows that \(d(jg^2)(z)=0\)
exactly when \(\mathfrak m(z)=0\), proving \eqref{eq:converse-kernel}.
\end{proof}

\begin{definition}[Linearization]\label{def:linearization}
For \(r=3,4,5\) and homogeneous \(x_i\in\bar H\), write
\[
  T_r(\mathbf x)=t[s_Hx_1|\cdots|s_Hx_r],
  \qquad
  e_r(\mathbf x)=\sum_{i=1}^r(r-i)|x_i|.
\]
Define the internal-degree-preserving linear map
\[
  \ell_r:\bigoplus_n\mathcal Q_{r-1}^{(n)}
  \longrightarrow C_r(H)
\]
on monomials with one factor in \(U_r\) by
\begin{equation}\label{eq:converse-linearization}
  \ell_r\bigl(P T_r(\mathbf x)\bigr)
  =
  (-1)^{(r-1)|P|+e_r(\mathbf x)}
  \rho(P)[x_1|\cdots|x_r],
  \qquad P\in\Lambda(U_1)\text{ homogeneous}.
\end{equation}
On monomials with at least two factors from
\(\bigoplus_{i\ge2}U_i\), set \(\ell_r=0\).
\end{definition}

Since \(U_i\) has resolution degree \(i-1\), a monomial of
resolution degree \(r-1\) containing exactly one factor from
\(\bigoplus_{i\ge2}U_i\) belongs to \(\Lambda(U_1)U_r\).
Thus the two prescriptions define \(\ell_r\) on its entire domain.
The sign \((-1)^{e_r(\mathbf x)}\) is the suspension sign of
Subsection \ref{subsec:obstructions}; it converts the shuffle
relations on \(s_H\bar H\) into the signed shuffle relations on
\(\bar H\).

\begin{remark}\label{rm:linear}
The lift \(\Gamma(z)\) lies in \(j(H)U_4\oplus U_2U_3\).
When \(b_3=0\), the evaluation defect satisfies
\[
  \mathfrak E|_{U_2U_3}=0.
\]
The motivation for this linearization is that we do not need to consider monomials lying in \(U_2U_3\) (see the proof of Proposition \ref{prop:converse}) in our comparison.
\end{remark}

\begin{lemma}\label{lem:converse-linearization}
Recall that \(\partial\) is the differential of the normalized Harrison chains. The following identities hold:
\[
 \ell_3d=-\partial_4\ell_4\quad\text{on }j(H)U_4\oplus U_2U_3,\qquad
 \ell_4d=-\partial_5\ell_5\quad\text{on }U_5 .
\]
\end{lemma}
\begin{proof}
For \(PT_r(\mathbf x)\), write \(p=|P|\) and
\(\theta=(r-1)p+e_r(\mathbf x)\). Using the formulas for \(D_1,D_2\)
in the proof of Proposition \ref{prop:eval-defect}, the terms of
\((-1)^\theta\ell_{r-1}d(PT_r)\) are as follows:
\[
\begin{array}{@{} c l l @{}}
\toprule
\text{Differential} & \text{Indices} & \text{Contribution} \\
\midrule
D_2 & (1,r-1)
    & -\rho(P)x_1[x_2|\cdots|x_r] \\[2pt]
D_1 & (i,i+1),\quad 1\le i<r
    & (-1)^{i+1}\rho(P)[x_1|\cdots|x_i x_{i+1}|\cdots|x_r] \\[2pt]
D_2 & (r-1,1)
    & (-1)^{r+1+|x_r|(p+\sum_{i<r}|x_i|)}
      x_r \rho(P)[x_1|\cdots|x_{r-1}] \\[2pt]
D_2 & (i,r-i),\quad 2\le i\le r-2
    & 0 \\
\bottomrule
\end{array}
\]
In the rows for \(D_2\), the pair \((i,r-i)\) means that the word
is cut after \(x_i\), giving the product
\(T_i(x_1,\ldots,x_i)T_{r-i}(x_{i+1},\ldots,x_r)\), with its
cobar sign. Here \(T_1(x)=j(x)\).
In the row for \(D_1\), the pair \((i,i+1)\) indicates that
the adjacent inputs \(x_i\) and \(x_{i+1}\) are replaced by
their product \(x_ix_{i+1}\) in \(H\).
The table proves both identities
on \(j(H)U_4\) and \(U_5\). For \(u\in U_2,\ w\in U_3\),
\(\ell_4(uw)=0\). Also \(\ell_3((du)w)=0\) since \(\rho(du)=0\),
and \(\ell_3(u\,dw)=0\) since \(dw\in j(H)U_2\).
\end{proof}

Since \(\ell_3Z=0\), Lemmas \ref{lem:finite-lift}
and \ref{lem:converse-linearization} define a map
\begin{equation}\label{eq:converse-psi}
 \Psi:\mathcal D^n(H)\longrightarrow H_4(C_\bullet(H)^n),\qquad
 z\longmapsto[\ell_4\Gamma(z)] .
\end{equation}
Indeed \(\partial_4\ell_4\Gamma(z)=-\ell_3 Z(z)=0\), and changing
a lift by \(dT\), \(T\in U_5\), changes its image by
\(-\partial_5\ell_5T\). Thus it is well defined.

\begin{proposition}\label{prop:converse-surjection}
The map \(\Psi\) is surjective in every internal degree.
\end{proposition}
\begin{proof}
Represent a cycle \(\zeta\in C_4(H)^n\) by
\(\sum_\nu c_\nu v_\nu[\mathbf x_\nu]\). Set
\[
 y=\sum_\nu c_\nu(-1)^{3|v_\nu|+e_4(\mathbf x_\nu)}
                 j(v_\nu)T_4(\mathbf x_\nu)\in j(H)U_4 .
\]
Then \(\ell_4 y=\zeta\). Choose a homogeneous basis
\(\{w_\lambda\}\) of \(U_3\). Since
\(dU_4\subset j(H)U_3+\mathcal G^2U_2\), there is a unique
decomposition
\begin{align*}
 dy &=\sum_\lambda p_\lambda w_\lambda+b,
 & b&\in j(H)\mathcal G^2U_2,\\
 p_\lambda &=c_\lambda1+j(a_\lambda)+j^2(q_\lambda),
 & a_\lambda&\in\bar H,\quad q_\lambda\in\mathcal G^2\bar H .
\end{align*}
Now we have \(\ell_3(dy)=-\partial_4\zeta=0\) and \(\ell_3(b)=0\).
In \(\ell_3(p_\lambda w_\lambda)\) the coefficient sign is
\((-1)^{2|p_\lambda|}=1\). The desuspended \(w_\lambda\) form a
basis of \(\bar H^{\otimes3}/\Sh_3\), so \(\rho(p_\lambda)=0\)
for each \(\lambda\). Thus \(c_\lambda=0\),
\(a_\lambda=-\mu(q_\lambda)\), and
\[
 p_\lambda=j^2(q_\lambda)-j(\mu(q_\lambda))=d g(e(q_\lambda)).
\]
Set \(y'=y-\sum_\lambda g(e(q_\lambda))w_\lambda\in j(H)U_4\oplus U_2U_3\).
Then \(\ell_4y'=\zeta\) and
\[
 dy'=b-\sum_\lambda(-1)^{|q_\lambda|-1}
             g(e(q_\lambda))\,dw_\lambda\in j(H)\mathcal G^2U_2,
\]
because \(dw_\lambda\in j(H)U_2\).
By Proposition \ref{prop:converse-kernel}, there exists \(z\in\mathcal D^n(H)\) such that \(dy'=Z(z)\).
Thus \(y'\) is a lift of \(Z(z)\) in \(j(H)U_4\oplus U_2U_3\), and the independence of
the lift in \eqref{eq:converse-psi} gives \(\Psi(z)=[\zeta]\).
\end{proof}

\subsection{Comparison with the fourth obstruction}

\begin{proposition}\label{prop:converse}
Under the hypotheses of Theorem \ref{thm:main},
\[
 \varepsilon\mathcal P_A(z)=-\overline\lambda_{m_4}(\Psi(z))
 \qquad (z\in\mathcal D^{N+2}(H)).
\]
In particular, if \(\varepsilon\mathcal P_A\) vanishes in degree
\(N+2\), then \(\kappa_4(A)=0\).
\end{proposition}
\begin{proof}
Write \(\Gamma(z)=\sum_\nu c_\nu j(v_\nu)T_4(\mathbf x_\nu)+\Gamma_{23}\)
with \(\Gamma_{23}\in U_2U_3\), allowing \(v_\nu=1\).
By \eqref{eq:forward-defect}, \(\mathfrak E(\Gamma_{23})=0\), and the Leibniz
rule for the defect gives
\[
 \mathfrak E\bigl(j(v)T_4(\mathbf x)\bigr)
 =(-1)^{|v|+e_4(\mathbf x)}\alpha(v)\alpha(m_4(\mathbf x)).
\]
It follows that
\[
\alpha\gamma^2(z)=\Phi Z(z)=D_A\Phi\Gamma(z)-\sum_\nu(-1)^{|v_{\nu}|+e_4(\mathbf x_{\nu})}c_{\nu}\alpha(v_{\nu})\alpha(m_4(\mathbf x_{\nu})).
\]
Taking cohomology classes and applying \(\varepsilon\), we obtain
\[
 \varepsilon\mathcal P_A(z)
 =-\sum_\nu c_\nu(-1)^{|v_\nu|+e_4(\mathbf x_\nu)}
            \varepsilon\bigl(v_\nu m_4(\mathbf x_\nu)\bigr)
 =-\lambda_{m_4}(\ell_4\Gamma(z)).
\]
The last equality uses \(3|v_\nu|\equiv |v_\nu|\pmod2\)
in \eqref{eq:converse-linearization}.

If \(\varepsilon\mathcal P_A(z)\) vanishes for every \(z\), surjectivity of \(\Psi\)
shows that \(\lambda_{m_4}\) vanishes on \(\ker\partial_4\) in
internal degree \(N+2\). Hence
\(L_0(\partial_4 c)=\lambda_{m_4}(c)\) defines a functional on
\(\im\partial_4\). Extend it linearly to \(C_3(H)^{N+2}\).
By Poincar\'e duality this extension equals \(\lambda_h\) for a
normalized Harrison map \(h:H^{\otimes3}\to H\) of degree \(-2\).
Equation \eqref{eq:converse-pd} with \(r=3\) now gives
\(\delta\widehat h=\widehat m_4=b_4\). Thus \(\kappa_4(A)=0\).
\end{proof}

Together with Proposition \ref{prop:forward} and the obstruction
criterion in Subsection \ref{subsec:obstructions}, this proves
Theorem \ref{thm:main}.

\section{Applications}\label{sec:applications}

\begin{proof}[Proof of Theorem \ref{cor:manifolds}]
Take a minimal Sullivan model \(A\) of \(M\) over \(\Q\).
Since \(M\) is connected, \(A^0=\Q\), and its cohomology
is \(H^*(M;\Q)\) with orientation
\(\varepsilon(x)=\langle x,[M]\rangle\) in degree \(N\).
Applying Theorem \ref{thm:main} gives Theorem \ref{cor:manifolds}.
\end{proof}

\begin{proof}[Proof of Theorem \ref{cor:spheres}]
Put \(R=H^*(B;\Q)\) and \(e_1=e(\pi_1)\).
Since the bundles are oriented, the fundamental groups of
their bases act trivially on the rational cohomology of
the fibres. The odd-sphere bundle model therefore applies
\cite[Theorem 20.3]{Halperin}.
Then we obtain a model of \(X\):
\[
  A=R\otimes\Lambda(u,v),\qquad
  D_A|_R=0,\qquad D_Au=e_1,\quad D_Av=e_2,
  \qquad |u|=p,\quad |v|=q.
\]
The total space \(X\) is closed, connected and oriented,
so \(H^*(A)\) satisfies Poincar\'e duality in dimension \(N\).

Give \(A\) the exterior length
\(A_j=R\otimes\Lambda^j(u,v)\), with \(A_j=0\) outside \(0\le j\le2\).
Since \(D_A A_j\subset A_{j-1}\), cohomology splits as
\[
 H=\bigoplus_{j=0}^2 H_j,\qquad
 H_j=\ker(D_A|_{A_j})/D_A A_{j+1}.
\]
Write \(A_j^k\) for the part of cohomological degree \(k\) and
exterior length \(j\). Since \(D_A:A_j^k\to A_{j-1}^{k+1}\), set
\[
  Z_j^k=\ker(D_A|_{A_j^k}),\qquad
  B_j^k=D_A(A_{j+1}^{k-1}).
\]
In each bidegree choose complements
\[
  Z_j^k=B_j^k\oplus I_j^k,\qquad
  A_j^k=Z_j^k\oplus L_j^k,
  \qquad 1\in I_0^0 .
\]
Thus \(I_j^k\) represents \(H_j^k\), and restriction gives an
isomorphism
\[
  D_A:L_{j+1}^{k-1}\xrightarrow{\;\cong\;}B_j^k .
\]
For \(b\in B_j^k\), \(a\in I_j^k\) and \(\ell\in L_j^k\), set
\[
  \alpha([a])=a,\qquad
  \pi(b+a+\ell)=[a],\qquad
  D_Ah(b+a+\ell)=b.
\]
Then we have
\[
  D_Ah+hD_A=\id_A-\alpha\pi,\qquad
  h:A_j^k\longrightarrow A_{j+1}^{k-1}.
\]

Apply the tree formula \eqref{eq:transfer-trees}.
The maps \(\alpha,\pi\) preserve exterior length, multiplication
adds the lengths of its inputs, and each of the \(r-2\) internal
edges raises exterior length by one. Hence
\begin{equation}\label{eq:sphere-length}
  m_r(H_{j_1},\ldots,H_{j_r})
  \subseteq H_{j_1+\cdots+j_r+r-2}
  \qquad (r\ge2).
\end{equation}
Since \(H_j=0\) for \(j>2\), it follows that \(m_r=0\) for \(r\ge5\).

Suppose condition \textup{(3)} holds. By Proposition \ref{prop:third},
choose a normalized cochain \(\xi_2\) with \(\delta\xi_2=b_3\).
A cochain has exterior-length shift \(a\) if it takes inputs of
lengths \(j_1,\ldots,j_r\) to length \(j_1+\cdots+j_r+a\).
Since \(b_2\) has shift zero, \(\delta=[\mathbf b_2,-]\) preserves
shift. Taking the shift-one part of \(\delta\xi_2=b_3\) therefore
gives a normalized cochain \(\xi_2\) of shift one.

Put \(X=\boldsymbol{\xi}_2\) and write
\(\operatorname{ad}_X(Y)=[X,Y]\). Conjugation by \(e^X\) gives
\[
  \mathbf b'=e^X\mathbf b e^{-X}
  =\sum_{n\ge0}\frac{1}{n!}\operatorname{ad}_X^n(\mathbf b).
\]
Each bracket with \(X\) raises arity by one, so the sum is finite
in each arity. Using
\(\operatorname{ad}_X(\mathbf b_2)=-\mathbf b_3\), we obtain
\[
\begin{aligned}
 \mathbf b'_3
   &=\mathbf b_3+\operatorname{ad}_X(\mathbf b_2)=0,\\
 \mathbf b'_4
   &=\mathbf b_4+\operatorname{ad}_X(\mathbf b_3)
       +\tfrac12\operatorname{ad}_X^2(\mathbf b_2)
     =\mathbf b_4+\tfrac12\operatorname{ad}_X(\mathbf b_3),\\
 \mathbf b'_5
   &=\mathbf b_5+\operatorname{ad}_X(\mathbf b_4)
       +\tfrac12\operatorname{ad}_X^2(\mathbf b_3)
       +\tfrac16\operatorname{ad}_X^3(\mathbf b_2)
\end{aligned}
\]
Every term in \(\mathbf b'_5\) has shift three and hence vanishes,
since \(H_j=0\) for \(j>2\).

Thus conjugation by \(e^X\) preserves \eqref{eq:sphere-length}
and shows that \(b'_4\) has shift two.

By Theorem \ref{thm:main}, choose a normalized \(\xi_3\) with
\(\delta\xi_3=b'_4\). Taking its shift-two part preserves this
equation. Put \(Y=\boldsymbol{\xi}_3\) and
\(\mathbf b''=e^Y\mathbf b'e^{-Y}\). Since \(Y\) has arity three,
\[
  \mathbf b''_3=\mathbf b'_3=0,\qquad
  \mathbf b''_4=\mathbf b'_4+[Y,\mathbf b_2]=0.
\]
The same arity and shift calculation gives \(b''_r=0\) for
\(r\ge5\). Thus the resulting structure has only \(b_2\).

When the \(C_\infty\) quasi-isomorphism to \(A\) is transported through both
conjugations, equation \eqref{eq:all-defects}
shows that its evaluation \(\Phi:\mathcal Q\to A\) is a CDGA map.
Since \(\rho\) is a quasi-isomorphism and
\([\Phi j(x)]=x\), both arrows in
\[
  (H,0)\xleftarrow{\ \rho\ }(\mathcal Q,d)
       \xrightarrow{\ \Phi\ }(A,D_A)
\]
are quasi-isomorphisms. Hence \(X\) is rationally formal.
This proves \textup{(3)}\(\Rightarrow\)\textup{(1)}.

Conversely, if \(X\) is rationally formal, the choice
\((\alpha,\gamma)=(\id,0)\) on \((H,0)\) gives
\(\mathcal F=\mathcal T=\mathcal P=0\).
Quasi-isomorphism invariance gives \textup{(2)}, and applying
\(\varepsilon\) in top degree gives \textup{(3)}.
\end{proof}

\subsection*{Higher obstructions and the construction of higher tensors}

Massey products suggest cochain formulas for higher tensors.
For example, on eight-dimensional \(E\)-manifolds Nagy and
Nordstr\"om identify the product of a fourfold Massey product
with a fifth class as an evaluation of the pentagonal tensor
\cite[Proposition 1.6]{NN}.
The choice of the relation space is essential: the
Bianchi-Massey tensor can be nonzero even when all ordinary
triple Massey products are trivial \cite[Section 1.5]{CN}.

Our proof suggests constructing these relation spaces directly
inside the bar-cobar model.
At arity four, Proposition \ref{prop:converse-kernel}
identifies the pentagonal relation space with the kernel of
\(d\) on \(j(H)\mathcal G^2U_2\).
In higher arity, one may similarly seek relations among the
closed elements in a chosen space of monomials determined by \(H\).

Suppose \(m_3,\ldots,m_{r-1}\) have been eliminated.
By Remark \ref{rem:higher-defect}, \(\Phi\) commutes with the
differential on \(\mathcal Q_{\le r-1}\), and its defect on
generators of length \(r\) represents \(m_r\).
One would seek a graded relation space \(\mathcal D_r(H)\)
and linear maps
\[
 Z_r:\mathcal D_r(H)\longrightarrow\mathcal Q_{\le r-1},
 \qquad
 \Gamma_r:\mathcal D_r(H)\longrightarrow\mathcal Q_{\le r},
 \qquad d\Gamma_r=Z_r,
\]
of degrees \(2-r\) and \(1-r\), respectively.
Then \(Z_r(z)\) and \(\Phi Z_r(z)\) are closed, and the
defect identity gives
\[
 [\Phi Z_r(z)]=-[\mathfrak E\Gamma_r(z)].
\]
It remains to prove that these classes are independent of
the choices of a strictly unital \(C_\infty\)-structure on \(H\)
and a compatible \(C_\infty\)-quasi-isomorphism to \(A\),
subject to \(m_3=\cdots=m_{r-1}=0\) and the fixed identification
\(H=H^*(A)\).
Once this is established, define
\[
 \mathcal P_r:\mathcal D_r^n(H)\longrightarrow H^{n+2-r}(A),
 \qquad
 \mathcal P_r(z)=[\Phi Z_r(z)].
\]
Eliminating \(m_r\) makes the defect vanish on
\(\mathcal Q_{\le r}\), so the displayed identity gives
\(\mathcal P_r=0\).

For the converse, suppose that \(H\) is a finite-dimensional
Poincar\'e duality algebra of formal dimension \(N\).
The corresponding graded Poincar\'e duality pairing identifies
the Harrison class of \(m_r\) with a functional
\[
 \overline{\lambda}_{m_r}:
 H_r\bigl(C_\bullet(H)^{N+r-2}\bigr)\longrightarrow\mathbb Q.
\]
Thus we focus on the induced functional on homology, without
requiring \(\lambda_{m_r}\) to vanish on all chains.
The next step would be to show that linearizing the lifts
defines a well-defined surjective map
\[
 \Psi_r:\mathcal D_r^{N+r-2}(H)
 \longrightarrow H_r\bigl(C_\bullet(H)^{N+r-2}\bigr),
\]
and to compare
\(\overline{\lambda}_{m_r}\circ\Psi_r\) with
\(\varepsilon\mathcal P_r\).
These two properties would give
\[
 \mathcal P_r=0
 \quad\Longrightarrow\quad
 \overline{\lambda}_{m_r}\circ\Psi_r=0
 \quad\Longrightarrow\quad
 \overline{\lambda}_{m_r}=0
 \quad\Longrightarrow\quad
 [m_r]=0.
\]
The last implication follows from Poincar\'e duality.
Finding relation spaces with the required lifts, choice
independence, and homological comparison is the remaining
problem in higher arity.

The exterior-length argument also provides a finite obstruction
range for certain models.
For a connected model
\[
 A=R\otimes\Lambda(u_1,\ldots,u_s),
 \qquad D_AR=0,\qquad D_Au_i\in R,
\]
where all \(u_i\) have positive odd degree, the argument in
the proof of Theorem \ref{cor:spheres} shows that \(A\)
is formal if and only if its successive obstruction classes
in arities \(3,\ldots,s+2\) vanish.

\end{document}